\documentclass[letterpaper]{article}

\usepackage{spconf}
\usepackage[dvipdfmx]{graphicx}
\usepackage{amsmath}
\usepackage{bm}
\usepackage{amsfonts}
\usepackage{amssymb}
\usepackage{multicol}

\renewcommand\Im{\operatorname{Im}}

\title{Upper-bounding bias errors in satellite navigation}
\name{Takashi Iwamoto}
\address{Advanced Technology R \& D Center\\
Mitsubishi Electric Corporation\\
Amagasaki, Hyogo 661-8661, Japan}
\begin{document}
\onecolumn
This manuscript is the accepted version of the paper:
\begin{quote}
T. Iwamoto, ``Upper-bounding bias errors in satellite navigation,'' 
to be presented at IEEE Workshop on Statistical Signal Processing, 
Gold Coast, Australia, 2014. 
\end{quote}

The following copyright notice applies:
\begin{quote}
\copyright 2014 IEEE. Personal use of this material is permitted. Permission from IEEE must be
obtained for all other uses, in any current or future media, including
reprinting/republishing this material for advertising or promotional purposes, creating
new collective works, for resale or redistribution to servers or lists, or reuse of any
copyrighted component of this work in other works.
\end{quote}
\vfill
\twocolumn
\maketitle
\begin{abstract}
A satellite navigation system for a safety critical application is
required to provide an integrity alert of any malfunction; the
probability that a navigation positioning error exceeds a given alert
limit without an integrity alert is required to be smaller than a given
integrity risk. So far, a little number of applications provide
integrity alerts, because signal propagation from a satellite to a
receiver depends on diversified phenomena and makes probabilistic
upper-bound of possible threats difficult. To widen application fields
of satellite navigation, two methods to upper-bound wide classes of bias
errors are shown in this paper. The worst bias error in a maximum
likelihood estimate caused by an interference signal within a given
small power is derived. A novel inequality condition with a clock bias
error and magnification coefficients that upper-bounds a horizontal
position error is presented. Robustness of the inequality condition is
numerically shown based on actual configurations of satellites.
\end{abstract}

\section{Introduction}
\label{sec:intro}

Satellite navigation systems have been tried on safety critical
applications in various fields. One of the most successful applications
is the Satellite Based Augmentation System (SBAS). It was commissioned
for aviation in USA in 2003~\cite{WAAS2008} and is planned to cover
large parts of the globe around 2020~\cite{Thomas2014}. SBAS provides
integrity alerts to contain the probability of Hazardously Misleading
Information below a given integrity risk. For example, the category of
Approach Procedures with Vertical Guidance I (APV-I) requires that the
probability that a vertical navigation positioning error exceeds the
alert limit of 50 m without integrity alerts should be less than $2
\times 10^{-7}$ per approach.

SBAS bases on over-bounding statistics that are designed to cover
potential threats with sufficient margins. For example, orbit and clock
errors of Global Positioning System (GPS) are investigated over a
considerable period and over-bounding distributions are
provided~\cite{Cohenour2011}. Multipath errors at antennas on airplanes
are also modeled and simulated~\cite{Lentmaier2008} and taken into
consideration of standard multipath powers as functions of parameters
such as an elevation angle of a satellite~\cite{ITU-R-P682}.

Besides aviation, applications to other fields including railways are
being tried~\cite{Ferrario2013}. Near surface, however, serious and
diversified errors caused by multipaths have been
reported~\cite{Steingass2004}. Since multipath errors show strong
dependences on each propagation environment, no distributions are
confirmed to over-bound multipath errors, so far. When sufficient
statistics are not available, additional sensors are expected to
overcome threats. It was pointed out that a precise clock of a receiver
works to detect an error in an estimate of clock bias that has
correlation with an error in an estimate of vertical
position~\cite{Bednarz2006}. In the same literature, a horizontal
position error was denoted as essentially uncorrelated with an estimate
of a clock bias error.

In this paper, two methods to upper-bound bias errors are shown. The
worst bias error in a maximum likelihood estimate caused by an
interference signal within a given small power is derived in a simple
expression in Section~\ref{sec:By power}. It gives fast estimation
of the worst bias error once the power of an interference signal becomes
available as in~\cite{ITU-R-P682}. In Section~\ref{sec:By geometry} a
novel inequality condition with a clock bias error and magnification
coefficients that upper-bounds a horizontal position error is
presented. Robustness of the inequality condition is numerically evaluated
based on an ephemeris of GPS satellites in Section~\ref{sec:Numerical
evaluation} and the conclusion follows.

\section{The worst bias error caused by an interference within a given
 small power}
\label{sec:By power}

A general framework of the worst bias error in a maximum likelihood
estimate caused by an interference signal within a given small power is
presented in~\cite{Wada2010}. For a cord spread signal as used in GPS,
it is further possible to express and evaluate the worst mode explicitly
as follows. Suppose that a signal is periodically sampled at $k T$ for
$k=1, 2,...$, where $T$ is a sampling period, and modeled by
\begin{equation}
 z(k T) = w(k T - \tau) + y(k T) + n(k T),
\end{equation}
where $w(k T - \tau) = r e^{i\phi} m(k T - \tau)$ is a known code $m(k
T - \tau) \in \mathbb{R}$ with an unknown delay $\tau \in \mathbb{R}$
multiplied by an amplitude $r$ and a phase $e^{i\phi}$, and $y(k T)$,
$n(k T) \in \mathbb{C}$ are respectively an interference signal and a
noise independently obeying a Gaussian distribution density
\begin{equation}
 p(n(k T)) = \frac{1}{\sqrt{2 \pi} \sigma}
 \exp \left( -\frac{|n(k T)|^2}{2\sigma^2}\right).
\end{equation}
The logarithmic likelihood function of parameters with respect to
samples is derived as $L_y = \sum_k \ell_{y,k}$ with
\begin{equation}
\resizebox{.89\hsize}{!}
{$\ell_{y,k} = 
  -\frac{|z(k T) - w(k T - \tau) - y(k T)|^2}{2 \sigma^2}
  -\frac{\ln (2 \pi \sigma^2)}{2}.
$}
\end{equation}
Let $\tau_0$ denote the maximum likelihood estimate of the delay for the
unperturbed ($y = 0$) system,
\begin{equation}
 \frac{\partial L_0}{\partial \tau}(\tau_0) = 0,
\end{equation}
and $\tau_0 + \delta\tau$ is the one for the perturbed system with a
small $\delta y$. The condition
\begin{equation}
 \frac{\partial L_{\delta y}}{\partial \tau}(\tau_0+\delta\tau) = 0
\end{equation}
gives an explicit representation
\begin{equation}
\resizebox{0.89\hsize}{!}
{$\begin{aligned}
  \sum_k &w'(k T - \tau_0-\delta\tau) \{z(k T) - w(k T -
  \tau_0-\delta\tau) -\delta y(k T)\}^{\ast}\\ & + \mbox{c.c.} = 0,
 \end{aligned}$}
\end{equation}
where $w'$ denotes the derivative of $w$, $z^{\ast}$ denotes the
complex conjugate of $z$, and c.c. denotes the complex conjugate of the
previous term.
Straightforward calculation up to the first order gives
\begin{align}
|\delta \tau| &= \left|\frac{\langle \delta y, w' \rangle }
 {\|w'\|^2+ \langle w-z , w''\rangle}\right|\\
 &\leq M_{\tau} \|\delta y\|,
\intertext{with the magnification coefficient defined by}
 M_{\tau} &:= \frac{\|w'\|} {\left| \|w'\|^2+\langle w-z, w'' \rangle\right|},
\end{align}
where $\langle a, b\rangle := \sum_k a(k)^{\ast} b(k)$ denotes the inner
product of complex vectors $a$ and $b$ of the same dimension, $\|\delta
y\| = \sqrt{\langle \delta y, \delta y\rangle}$ is the square root of a
power of an interference signal $\delta y$, and $w''$ denotes the second
derivative of $w$. The equality holds for the interference signal
$\delta y$ parallel to $w'$. This worst mode corresponds to the
Goldstone mode with respect to the time translational
symmetry~\cite{Anderson1984} and gives fast estimate of the worst error
caused by an interference within a given small power regardless of
details, once the power of interference becomes available.

\section{Upper-Bounding bias errors in horizontal position estimates}
\label{sec:By geometry}

The observation equation of a pseudo range
\begin{equation}
 \delta\rho_j = \langle\bm{g}_j, \delta \bm{x}\rangle + \delta b + \epsilon_j,
 \label{eq:observation}
\end{equation}
is derived by linearizing a pseudo range $\rho_j$ from the $j$-th
navigation satellite ($j=1, 2,\ldots$) reaching a base point $\bm{x}_0$
at a base time $t_0$ for a neighborhood $(\bm{x}_0+\delta \bm{x},
t_0+\delta b/c) \in \mathbb{R}^4$, where a geometry vector $\bm{g}_j \in
\mathbb{R}^3 $ is the unit direction vector to the $j$-th satellite form
the base point $\bm{x}_0$ multiplied by $-1$, $\epsilon_j$ is an
correction term on the pseudo range form the $j$-th satellite, $b$ is
called a bias, and $c$ denotes the light velocity.

Suppose that a position of an antenna is limited in a curve $\bm{x}(s) =
(x(s), y(s), z(s)) \in \mathbb{R}^3$ that is continuous and an almost
everywhere differentiable function of the distance $s$ from the base
point (Fig.~1 (a)). Let $\bm{U}$ and $\bm{V}$ be the tangential unit
vector of the curve and the unit vector to the center of the osculating
circle that is tangential to the curve at the point $\bm{x}(s)$
respectively, and $\bm{W}$ be the unit vector that consists of a
orthogonal frame with $\bm{U}$ and $\bm{V}$, which is called the Frenet
frame (Fig.~1 (b)).
\begin{figure}[t]
\begin{minipage}{\linewidth}
  \centering
  \centerline{\includegraphics[width=0.8\linewidth]{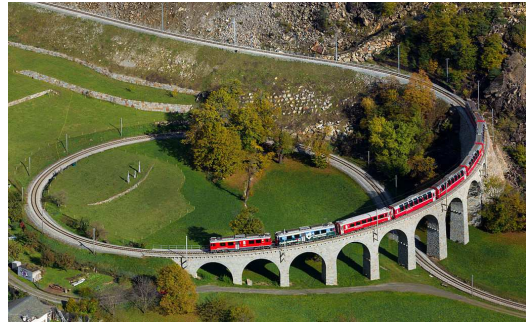}}
  \centerline{(a) An example of a railway curve (Brusio spiral viaduct~\cite{Gubler}).}\medskip
\end{minipage}\\
\begin{minipage}{\linewidth}
  \centering
  \centerline{\includegraphics[width=0.4\linewidth]{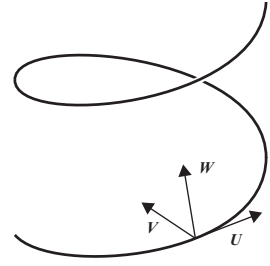}}
 \centerline{(b) The Frenet frame on a curve.}\medskip
\end{minipage}
\caption{An example
of a frame on a railway curve.}  \label{fig:the Frenet frame}
\end{figure} 
Against this frame a deviation is expressed as $\delta \bm{x} = \delta u
\bm{U} + \delta v \bm{V} + \delta w \bm{W}$.  The curve is approximated
by the osculating circle, whose radius is denoted by $R$, in a
neighborhood of $\bm{x}(s)$ as
\begin{equation}
 (u(s), v(s)) = R \left( \sin \frac{s}{R}, 1-\cos \frac{s}{R}\right).
 \label{eq:approximation}
\end{equation}
Under this approximation, a projection from the point $(u,v)$ measured
from the navigation satellites to the distance parameter $s = R \arctan
\left(u/(R - v)\right)$ is derived for the domain $v < R$
and linearized for small deviations
\begin{equation}
\resizebox{0.89\hsize}{!}
{$\delta s = \frac{R (R-v) (u^2 + (R - v)^2)}{(R - v)^4} \delta u +
 \frac{R u (u^2 + (R - v)^2)} {(R - v)^4}\delta v.
$}
\end{equation}

When three satellites $j=1, 2, 3$ are available, it is derived from 
the observation equation 
\begin{equation}
 \left[
 \begin{array}{c}
  \delta \rho_1\\ \delta \rho_2\\ \delta \rho_3
 \end{array}
 \right]
 = \left[
 \begin{array}{ccc}
  f_1 & h_1 & 1\\
  f_2 & h_2 & 1\\
  f_3 & h_3 & 1\\
 \end{array}
 \right]
 \left[
 \begin{array}{c}
  \delta u\\ \delta v\\ \delta b
 \end{array}
 \right]
 +
 \left[ 
 \begin{array}{c}
  \epsilon_1\\ \epsilon_2\\ \epsilon_3
 \end{array}
 \right],
\end{equation}
where directional cosines are denoted by $f_j :=
\langle\bm{g}_j, \bm{U}\rangle$, $h_j := \langle \bm{g}_j \bm{V}\rangle$
respectively. 
When the condition on satellites configuration
\begin{equation}
 D := f_1 h_2 - f_2 h_1 + f_2 h_3 - f_3 h_2 + f_3 h_1 - f_1 h_3 \ne 0
\end{equation}
holds, the equation
\begin{equation}
\resizebox{0.89\hsize}{!}
{$\arraycolsep=1pt %
\left[
  \begin{array}{c}
   \delta u\\ \delta v\\ \delta b
  \end{array}
  \right]\negthickspace%
 =\negthickspace\frac{1}{D}\negthickspace%
  \left[
  \begin{array}{ccc}
  h_2 - h_3 & h_3 - h_1 & h_1 - h_2 \\
  f_3 - f_2 & f_1 - f_3 & f_2 - f_1 \\
  f_2 h_3 - f_3 h_2 & f_3 h_1 - f_1 h_3 & f_1 h_2 - f_2 h_1
 \end{array}
 \right]\negthickspace\negthickspace\left[
 \begin{array}{c}
  \delta \rho_1 - \epsilon_1\\ 
  \delta \rho_2 - \epsilon_2\\ 
  \delta \rho_3 - \epsilon_3\\ 
  \end{array}
  \right]$}
\end{equation}
holds. When the correct value of each correction $\epsilon_j$ is
substituted, this equation represents the change of coordinates
$(\delta u, \delta v, \delta b)$ against the change of each pseudo range
$\delta \rho_j$. If some correction is not included correctly, then the
output contains error. 

Empirical distributions of orbit and clock errors of navigation
satellites shows good convergence to the stationary
one~\cite{Cohenour2011}. There are also augmentation services to monitor
orbit and clock errors. Thus it is well assumed that those errors are
over-bounded by a certain distribution and its scale can be neglected 
later in this paper, when larger residuals are discussed.

Besides above mentioned, ionospheric delay, tropospheric delay,
diffraction of signals, and reflection of signals are major sources of
errors. All these physical effects make positive contribution to the
correction term $\epsilon_j > 0$; if not all are substituted, there
remains error $r_j := \delta \rho_j - \epsilon_j > 0$.

For $j=1, 2, 3$ complex numbers $z_j := f_j + i h_j$ and $z_4 = z_1$
satisfy an equation $\Im(z_j^\ast z_{j+1}) = f_j h_{j+1} - f_{j+1}
h_{j}$. There exists complex numbers that satisfy the condition $\Im
(z_1^\ast z_1 z_2^\ast z_2 z_3^\ast z_3) = 0$ subject to the condition
$\Im(z_j^\ast z_{j+1}) > 0$ under suitable exchanges of suffixes. For
satellites satisfying the above condition, an upper bounding inequality
for an error along a track 
\begin{align}
 |\delta u| &= 
  \left|
   \frac{(h_2 - h_3) r_1 + (h_3 - h_1) r_2 + (h_1 - h_2) r_3}{D}
  \right| \\ 
  &\leq M_u  |\delta b| 
\intertext{and one for perpendicular to a track}
 |\delta v| &\leq M_v  |\delta b|
\end{align}
are derived, where magnification coefficients are defined by
\begin{equation}
\resizebox{0.89\hsize}{!}
{$\begin{aligned}
 M_u &:= \frac{\max(|h_2 - h_3|, |h_3 - h_1|, |h_1 - h_2|)} {\min(|f_2 h_3
            - f_3 h_2|, |f_3 h_1 - f_1 h_3 |, | f_1 h_2 - f_2 h_1|)}, \\
 M_v &:= \frac{\max(|f_2 - f_3|, |f_3 - f_1|, |f_1 - f_2|)} {\min(|f_2 h_3
            - f_3 h_2|, |f_3 h_1 - f_1 h_3 |, | f_1 h_2 - f_2 h_1|)}.
\end{aligned}$}
\end{equation}

On the curve whose radius $R$ is large enough, positioning error in the
direction perpendicular to the curve does not affect an along track error
effectively. In those cases, by introducing a virtual satellite $j=3$ to
cancel the perpendicular error, only two physical satellites are needed
for positioning and its evaluation as shown below.  For a sufficient
large $h_3$, the equation
\begin{equation}
 \left[
 \begin{array}{c}
  \delta \rho_1\\ \delta \rho_2\\ \delta \rho_3
 \end{array}
 \right]
 = \left[
 \begin{array}{ccc}
  f_1 & h_1 & 1\\
  f_2 & h_2 & 1\\
  0 & h_3 & 0\\
 \end{array}
 \right]
 \left[
 \begin{array}{c}
  \delta u\\ \delta v\\ \delta b
 \end{array}
 \right]
 +
 \left[ 
 \begin{array}{c}
  \epsilon_1\\ \epsilon_2\\ \epsilon_3
 \end{array}
 \right]
\end{equation}
is considered. 
For the configuration of satellites satisfying the determinant condition
$D' := f_2 - f_1 \ne 0$, the equation 
\begin{equation}
\resizebox{0.89\hsize}{!}
{$\arraycolsep=1pt %
 \left[
 \begin{array}{c}
  \delta u\\ \delta v\\ \delta b
 \end{array}
 \right]
 = \frac{1}{D'} \left[
 \begin{array}{ccc}
  - 1 & 1 & (h_1 - h_2)/h_3 \\
  0 & 0 & (f_2 - f_1) h_3 \\
  f_2 & - f_1 & (f_1 h_2 - f_2 h_1)/h_3
 \end{array}
 \right]\negthickspace\negthickspace\left[
 \begin{array}{c}
  \delta \rho_1 - \epsilon_1\\ 
  \delta \rho_2 - \epsilon_2\\ 
  \delta \rho_3 - \epsilon_3\\ 
 \end{array}
 \right]$}
\end{equation}
holds. Arguments parallel to the previous on gives necessary conditions,
either $f_2 > 0 $ and $f_1 < 0$, or $f_2 < 0 $ and $f_1 > 0$. When one
of the conditions is holds, formally in the limit of $h_3 \rightarrow \infty$
with a magnification coefficient
\begin{align}
 M_s & := \frac{1}{\min(|f_1|, |f_2|)}, 
\intertext{an upper-bounding inequality along a track is derived as} 
|\delta s| &\leq M_s |\delta b|.
\end{align}
Note that this result is
consistent with one from a system
\begin{equation}
  \left[ 
 \begin{array}{c}
  \delta \rho_1\\ \delta \rho_2
 \end{array}
 \right]
 = \left[
 \begin{array}{cc}
  f_1 & 1\\
  f_2 & 1
 \end{array}
 \right]
 \left[
 \begin{array}{c}
  \delta s\\ \delta b
 \end{array}
 \right]
 +
 \left[ 
 \begin{array}{c}
  \epsilon_1\\ \epsilon_2
 \end{array}
 \right].
\end{equation}

\section{Numerical evaluation}
\label{sec:Numerical evaluation} 

\begin{figure}[t]
  \centering \centerline{\includegraphics[width=0.85\linewidth]
  {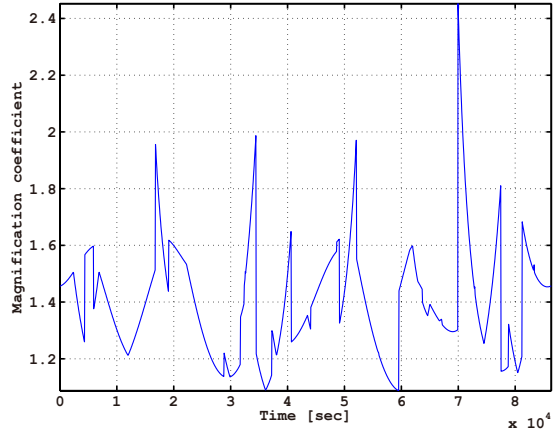}}
  \caption{\label{fig:MagnificationTime} The magnification coefficient
  calculated from the ephemeris on 25 July 2013 against every 60
  seconds.}
\end{figure}

\begin{figure}[t]
  \centering \centerline{\includegraphics[width=0.85\linewidth]
  {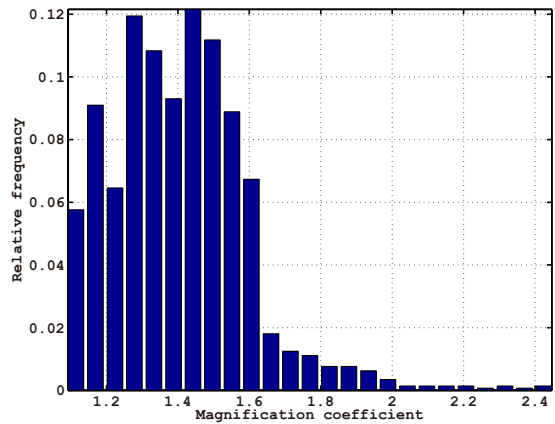}}
  \caption{\label{fig:MagnificationBar} Relative frequency of the
  magnification coefficients calculated for the ephemeris on 25 July
  2013.}
\end{figure}

The magnification coefficient $M_s$ defined in the previous section is
calculated from an ephemeris on 25 July 2013 (Universal Time). Satellite
positions are converted into local coordinates at the point on a railway
(latitude 34.75337, longitude 135.42783, height 3.7m) laid almost
straight in East-West direction near our factory in Amagasaki,
Japan. Magnification coefficients at every 60 seconds are calculated
on satellites with an elevation mask of 15 degree, and plotted
against time in Fig. \ref{fig:MagnificationTime}.
Sharp transitions in magnification coefficient values are caused by
appearance or disappearance of satellites and thus depend on the
selection of an elevation mask. Relative frequency of the magnification
coefficients are plotted in Fig. \ref{fig:MagnificationBar}.
It is observed that most magnification coefficients are contained under
the value 2 and concentrated under the value 1.6.

The accuracy of a local clock depends on its mechanism and calibration
scheme. There are wide range of choices including a rubidium clock whose
frequency stability is demonstrated less than $4 \times 10^{-12}
\tau^{-1/2}$ up to $10^{4}$ sec~\cite{Affolderbach2006}. 

\section{CONCLUSION}

In this paper two methods to upper-bound bias errors have been
shown. The worst bias error in maximum likelihood estimates caused by an
interference signal within a given small power is represented. It gives
a fast evaluation of an error caused by an interference regardless of
details, once the power of interference becomes available. A novel
inequality condition that upper-bounds a horizontal position error with
a clock bias error and magnification coefficients is also
shown. Robustness of the inequality conditions are shown based on actual
constellation of satellites.  Although optimization of a local clock as
well as its calibration scheme and evaluation of system performance are
remained for future work, these upper-bounds are expected to provide
concrete bases for wide fields of safety critical applications.


\begin{thebibliography}{99}
 \bibitem{WAAS2008} \textit{Global positioning system with wide are
	 augmentation system (WAAS) performance standard}, Federal
	 Aviation Administration, Washington D.~C., 2008.
 \bibitem{Thomas2014} D.~Thomas, ``GNSS Evolution Program Update,'' in
	 \textit{Interoperability working group meeting (IWG)}, 2014.
 \bibitem{Cohenour2011} C.~Cohenour and F.~van Graas, ``GPS Orbit and
	 Clock Error Distributions,'' \textit{NAVIGATION: Journal of the
	 Institute of Navigation}, vol.~58, no.~1, pp.~17--28, 2011.
 \bibitem{Lentmaier2008} M.~Lentmaier, B.~Krach, T.~Jost, A.~Lehner, and
	 A.~Steingass, \textit{Assessment of multipath in aeronautical
	 environments,} Citeseer, 2008.
 \bibitem{ITU-R-P682} \textit{Propagation data required for the design
	 of Earth-space aeronautical mobile telecommunication systems},
	 International Telecommunication Union, ITU-R P.~682-3, 2012.
 \bibitem{Ferrario2013} A.~Ferrario, L.~Marradi, P.~Iacone, and
	 A.~Galimberti, ``Multi-constellation GNSS Receiver for Rail
	 Applications,'' in \textit{ION GNSS+}, 2013.
 \bibitem{Steingass2004} A.~Steingass and A.~Lehner, ``Measuring the
	 navigation multipath channel a statistical analysis,'' in
	 \textit{ION GPS 2004 Conference}, 2004.
 \bibitem{Bednarz2006} S.~Bednarz and P.~Misra, ``Receiver clock-based
	 integrity monitoring for GPS precision approaches,''
	 \textit{IEEE Transactions on Aerospace and Electronic
	 Systems}, vol.~42, no.~2, pp.~636--643, 2006.
 \bibitem{Wada2010} T.~Wada and T.~Iwamoto, ``Evaluation of bias errors
	 in positioning a radio transmitter with a lot of interference
	 waves,'' in \textit{Phased Array Systems and Technology
	 (ARRAY)}, IEEE, 2010, pp.~1033--1038.
 \bibitem{Anderson1984} P.~W.~Anderson, \textit{Basic notions of condensed
	 matter physics}, Benjamin, 1984.
 \bibitem{Gubler} Gubler, ``Picture taken near Brusio, Switzerland,''
	 http://www.bahnbilder.ch/picture/11543, Creative Commons
	 Attribution 3.0 Unported. 
 \bibitem{Affolderbach2006} C.~Affolderbach, F.~Droz, and G.~Mileti,
	 ``Experimental Demonstration of a Compact and High-Performance
	 Laser-Pumped Rubidium Gas Cell Atomic Frequency Standard,''
	 \textit{IEEE Transactions on Instrumentation and Measurement},
	 vol.~55, no.~2, pp.~429--435, 2006.
\end{thebibliography}
\end{document}